\documentclass[10pt]{amsart}
\begin{document}
\title{Rational singular loci of nilpotent varieties}
\author{William M. McGovern}
\begin{abstract}
We present two methods for computing the rational singular locus of the closure of a nilpotent orbit in a complex semisimple Lie algebra and give a number of interesting examples.
\end{abstract}
\maketitle

\section{Introduction}
Let $G$ be a complex reductive group with Lie algebra $\frak g$.  Let $\mathcal N$ be the nilcone of $\frak g$ and let $X=\overline{G\cdot e}$ be a nilpotent variety in $\mathcal N$ (the closure of a nilpotent orbit $\mathcal O=G\cdot e$).  Recall that a point $x$ of $X$ is said to be rationally smooth if for all $y$ in a neighborhood of $x$ in the complex topology we have $H_y^m(X) = 0$ for all $m\ne 2\dim X$ and $H_y^{2\dim X}(X) = \mathbb Q$, where $H_y$ denotes cohomology with support in $\{y\}$; here we can use either singular or intersection cohomology \cite{GM83}.    The rational singular locus Rat Sing~$X$ consists of all points of $X$ at which $X$ is rationally singular (not smooth); this is a closed subset of $X$.  Rational smoothness has played a major role in representation theory ever since the seminal work of Kazhdan and Lusztig in the seventies \cite{KL79}; it also provided one of the first applications of intersection cohomology.  We will describe two methods of computing Rat Sing~$X$ and give some examples.
  
\section{First method--Brion}

Our first method applies to very general varieties $X$; they need not admit dense orbits under the action of any reductive group, but they must carry an action of a torus $T$ of dimension at least two.    This method uses techniques of Brion first developed in \cite{Br99}.  

\newtheorem*{theorem}{Theorem}
\begin{theorem}
Suppose that the variety $X$ admits an action of a torus $T$ such that $x\in X$ is an attractive fixed point of the $T$ action (so that all weights of $T$ on the tangent space $T_x X$ lie in an open half-space.  Then $X$ is rationally smooth at $x$ if and only if 

\begin{itemize}
\item Some punctured neighborhood of $x$ in $X$ is rationally smooth.

\item $X^{T'}$ is rationally smooth at $x$ for every subtorus $T'$ of codimension one in $T$.

\item $\dim_x X = \sum_{T'} \dim_x X^{T'}$ (sum over all subtori $T'$ of codimension one in $T$)

\end{itemize}
In addition, for any $T$-fixed point $x$ (attractive or not), the second condition is necessary for rational smoothness at $x$ and the third one becomes a necessary condition if $\dim_x X^T$ is subtracted from the left side and from every term in the right side.
\end{theorem}

For the proof see \cite[1.1,1.4]{Br99}.  Here the first of these conditions is typically hard to verify, but essential:  the second and third (in the absence of the first) might hold for one torus $T$ but fail for another, and in addition the first condition is needed to check that the cohomology vanishes at every point in a neighborhood of the one in question.  We can apply these conditions in particular to a nilpotent variety $\overline{G\cdot e}$ and a point $x$ on it, provided that a torus $T$ of dimension at least two fixes $x$:  this will occur whenever $x$ lies in the derived subalgebra of a Levi subalgebra of $\frak g$ of corank at least two.  The case $x=0$ is especially easy; here we may take the torus to be the direct product of a maximal torus $T$ of $G$ and a 1-torus $\mathbb C^*$, the latter acting on $\frak g$ by scalar multiplication.  Then $0$ is always an attractive fixed point for the action of this larger torus.

\newtheorem*{corollary}{Corollary}
\begin{corollary}
The full nilcone $\mathcal N$ is rationally smooth.  Every nilpotent variety $\bar{\mathcal O}$ different from $\mathcal N$ and 0 is rationally singular at 0, with two exceptions, namely the closures of the minimal orbit in types $C_n$ and $G_2$.
\end{corollary}

The rational singularity result is a straightforward consequence of the last of Brion's criteria for rational smoothness:  unless $\mathcal O$ is the minimal orbit, the subtori $T'$ contributing to the right side are exactly the centralizers of the positive root spaces in $\frak g$ and every nonzero term on the right side equals 2, whence the right side equals the dimension of $\mathcal N$ rather than $\bar{\mathcal O}$.  Even if $\mathcal O$ is minimal, the right side winds up being too large, except in types $C_n$ and $G_2$.  In these cases all three of Brion's criteria hold, the last two by direct calculation and the first one since there is only one singular point.  Thus these varieties are rationally smooth.  That $\mathcal N$ is rationally smooth is an old result of Borho and Macpherson \cite{BM83}; we will sketch their argument in the next seciton.  This result can also be proved using Brion's techniques, by constructing a slice of $\mathcal N$ in the sense of \cite[2.1]{Br99}.  This can be done inductively, starting with the full nilcone $\mathcal N'$ for $\frak{sl}_2$, which is well known to be rationally smooth.

It seems difficult in general to prove rational smoothness of nilpotent varieties using these techniques (because of the difficulty of verifying the first of Brion's conditions), but one can for example show that the rational singular locus of any spherical nilpotent variety in type $C$ lies one step below its boundary (in the chain of spherical varieties, ordered by inclusion; recall that a nilpotent orbit or its variety is called spherical if it admits a dense suborbit under the action of a Borel subgroup).  The same result holds for the largest spherical variety in type $A$ in odd rank; for other spherical varieties in this type, the rational singular locus coincides with the boundary.

\section{Second method--Borho-MacPherson}

We turn now to the second method, which is due to Borho and MacPherson and applies specifically to nilpotent varieties \cite{BM83}.  It uses the heavier machinery of intersection homology, but is able to compute the dimensions of the cohomology groups (rather than just determining whether they vanish outside the top degree).  We begin by invoking the Springer correspondence:  given $e$ the Springer fiber $\mathcal B_e$ of Borel subalgebras containing $e$ is such that its (singular) cohomology groups carry commuting actions of the component group $A(G\cdot e)$ of the centralizer of $e$ in $G$ and the Weyl group $W$.  The representation $\sigma_e$ of $W$ on the $A(G\cdot e)$-fixed vectors in the top cohomology group is then irreducible \cite{S78}.  Given now another nilpotent element $x$, the Borho-MacPherson Decomposition Theorem \cite{BM83} then implies that
$$
\dim IH^i_x(\bar{\mathcal O},\mathbb Q) = \dim\hom (\sigma_e,H^i(\mathcal B_x.\mathbb Q))
$$
\noindent where as usual $\bar{\mathcal O}$ denotes the closure of the orbit through $e$.  Here the left hand side denotes the intersection homology groups of $\bar{\mathcal O}$, which are naturally isomorphic to its cohomology groups in the complementary degree.  Hence $\bar{\mathcal O}$ is rationallly smooth at $x$ if and only if $IH^i_x(\bar{\mathcal O},\mathbb Q)$ vanishes in positive degree.  If $\bar{\mathcal O} =\mathcal N$, then $\sigma_e$ is trivial and an old computation of Lusztig shows that $\sigma_e$ occurs once in $H^i(\mathcal B_x,\mathbb Q)$ if $i=0$ and not at all if $i>0$ (for any $x$), so $\mathcal N$ is rationally smooth.   In general it is not easy to compute the module structure of $H^*(\mathcal B_x,\mathbb Q)$, but explicit tables are available in the exceptional cases \cite{BS84} while in the classical cases one has algorithms due to Shoji and Lusztig \cite{L81,FMM13,Sh83}.  An easy special case occurs when $x$ is regular nilpotent in some proper Levi subalgebra $\frak l$ of $\frak g$:  if $W_L\subset W$ is the Weyl group of $\frak l$, then $H^*(\mathbb B_x,\mathbb Q)$ is just the permutation representation of $W$ on $W/W_L$ \cite{AL82}.  This permutation representation can be computed by the Littlewood-Richardson rule in the classical cases and tables of Alvis in the exceptional ones.

To determine rational smoothness at $x$ it remains to compute $IH^*_y$ for all $y$ in a neighborhood of $x$; by $G$-equivariance it suffices to compute this group for each of the finitely many orbits $G\cdot y$ whose closures lie between $\overline{G\cdot x}$ and $\bar{\mathcal O}$.  One could ask whether cohomology vanishing at $x$ implies cohomology vanishing in a neighborhood (as it does for Schubert varieties).  The answer is no, already in type $C_3$:  there one computes that the cohomology of $\bar{\mathcal O}_{3,3}$ vanishes at points of $\mathcal O_{2,1^4}$, but not at points of $\mathcal O_{2^2,1^2}$, where $\mathcal O_{\lambda}$ denotes the orbit with partition $\lambda$ and exponents in partitions as usual denote repeated parts.  It also fails in type $D_4$:  the cohomology of $\bar{\mathcal O}_{5,3}$ vanishes at points of $\mathcal O_{3,2^2,1}$ but not at points of $\mathcal O_{3^2,1^2}$.  But it holds in type $A$, as follows from the combinatorics of Kostka numbers and the fact that every nilpotent element in that type is regular in some Levi subalgebra.

We conclude with an example of a rational singular locus of codimension 2 (this cannot happen for Schubert varieties, or for closures of $K$-orbits in the flag variety $G/B$, where $K$ is a symmetric subgroup of $G$).  Take $\mathcal O$ to be the orbit with Bala-Carter label $A_4$ in type $E_6$; this has dimension 60.  Applying the second method and using \cite{BS84} we compute that the rational singular locus of $\bar{\mathcal O}$ coincides with its boundary and has dimension 58; this is the closure of the orbit $\mathcal O'$ with Bala-Carter label $D_4(a_1)$.  In this case $\sigma_e$ occurs with multiplicity 3 in $H^*(\mathcal B_x)$ (where $e\in\mathcal O,x\in\mathcal O'$).  In fact $\sigma_e $ is paired with the two-dimensional reflection representation of the component group $A(\mathcal O')$ in the Springer correspondence.  The corresponding phenomenon also occurs for one the 42-dimensional orbits in type $F_4$ and the (unique) 40-dimensional orbit contained in its closure; there the component group of the smaller orbit is the symmetric group $S_4$ and so once again the multiplicity of the relevant Springer representation is larger than one.

\end{document}